\documentclass[12pt,reqno]{amsart}

\newtheorem{theorem}{Theorem}%[section]
\newtheorem{corollary}{Corollary}%[section]

\newcommand{\ZZ}{{\mathbb{Z}}}

\newcommand{\RR}{{\mathbb{R}}}

\newcommand{\NN}{{\mathbb{N}}}
\newcommand{\LL}{{\mathcal{L}}}

\newcommand{\G}{\mathcal{G}}
\newcommand{\J}{\mathcal{J}}
\newcommand{\I}{\mathcal{I}}

\newcommand{\D}{{\mathcal{D}}}
\newcommand{\M}{\mathcal{M}}

\newcommand{\norm}[1]{\left|\hspace*{2.5pt} \!\!\left|
#1\right|\hspace*{2.5pt} \!\!\right|} 

\newcommand{\card}{\operatorname{card}}

\numberwithin{equation}{section}
\numberwithin{figure}{section}

\begin{document}
\author{Cristian Cobeli
}

 \address{Cristian Cobeli,
 Institute of Mathematics of the Romanian Academy,
 P. O. Box \mbox{1-764}, Bucharest 70700,
 Romania.}
 \email{cristian.cobeli@imar.ro}

\title[On the discrete logarithm problem]
{On the discrete logarithm problem}

\subjclass[2000]{
Primary
11A07; %Congruences; primitive roots; residue systems
Secondary
11B50, %Sequences (mod m)
%11B57 Farey sequences; the sequences 1
11L07 %Estimates on exponential sums 
%11T23 %Exponential sums 
}
\thanks{Key Words and Phrases: Discrete logarithms, exponential sums,
  characters, primitive roots}

%%%%%%%%%%%%%%%%%%%%%%%%%%%%%%%%%%%%%%%%%% %%%%%%%%%
%\date{\today}

\begin{abstract}
Let $p>2$ be prime and $g$ a primitive root modulo $p$.
We present an argument for the fact that discrete logarithms of the
numbers in any arithmetic progression are uniformly distributed in
$[1,p]$ and raise some questions on the subject.
\end{abstract}

\maketitle
%%%%%%%%%%%%%%%%%%%%%%%%%%%%%%%%%%%%%%%%%% 
\section{Introduction}\label{sec:Introduction}
Before the middle of the last century, discrete logarithms were just
common tools used to perform calculations in finite fields. Then, with
the development of cryptography, their importance raised considerably,
especially after Diffie and Hellman~\cite{DH} created the key exchange
algorithm, the first practical public key cryptosystem. Many cryptosystems, such
as the Diffie-Hellman key agreement and its derivatives, ElGamal public-key
encryption, ElGamal signature scheme and its variants, DSA, etc. (see
\cite{ElGamal}, \cite{DSA}, \cite{Menzes}) are based on the assumption that
discrete logarithms are hard to compute. Considerable efforts have been made to
find algorithms that speed up the calculation of discrete logarithms, but nobody
knows how one could prove that a very fast algorithm does not exist. 

A strong argument would require proofs for the random distribution
characteristic of the set containing the discrete logarithms of the elements of
a ``regular'' subset of $[0,p-1]$ (a subinterval being just the first try),
when $p\to\infty$. This feature is suggested by numerical evidences for small
$p$ and by most of the work done around the cryptosystems based on discrete
logarithms (see \cite{Menzes} and \cite{Odlyzko} and the references
within). Recently, Banks and Shparlinski \cite{BS} have obtained nice results in this direction.

Discrete logarithms can be defined in general groups, but we reduce
here only to the group $\G=\ZZ/p\ZZ$ of residue classes modulo a prime
$p>2$. 
Given any $g\in \G$ and $n\in\NN$, let $g^n:=\underbrace{g\cdots g}_{n\ g's}$ be
the discrete exponentiation function. We will assume that $g$ is a
generator of $\G$, that is, $g$ is a primitive root modulo $p$. Then,
for any $x\in \G$, \emph{the discrete logarithm problem}
requires to find the smallest integer with the property that
$g^n\equiv x \pmod p$. Since $g$ is a primitive root, the power $n$
always exists in the interval $[0,p-1]$. We denote it by $n=\log_g
x=\log x$, and call it the \emph{the discrete logarithm} of $x$ to
base $g$.

Notice that the discrete logarithm function is the inverse of the
discrete exponentiation function and it has the properties $\log 1 = 0$
and $\log xy=\log x+\log y \pmod{p-1}$, for any $x,y\in \G$.

Let $a\ge 0$, $r>0$, $N>0$ be integers, and set 
$\J=\big\{a+r, \dots, a+Nr\big\}\subset[1,p-1]$. Denote
  \begin{equation*}
    \LL(g,\J)=\LL(g):=\big\{\log_g (a+jr)\colon\ 1\le j\le N\big\}
  \end{equation*}
and
  \begin{equation*}
    \M(g,\J)=\M(g):=\Big\{\frac{\log_g (a+jr)}{p-1}\colon\ 1\le j\le N\Big\}\,,
  \end{equation*}
the image of $\J$ in the torus $\RR/\ZZ$.
Then, any property regarding the spreading of the elements of $\LL(g,\J)$ over
$[0,p-1]$ transfers into a similar one regarding the elements of $\M(g,\J)$
over the torus, and
conversely. Since our aim is to understand what happens when $p$ gets large,
and it is more convenient to work within the bounded space $\RR/\ZZ$, in the
following our focus will concentrate mainly on $\M(g,\J)$.

The \emph{discrepancy} of $\M(g)$ is defined by
  \begin{equation*}
    \D\big(\M(g); \alpha,\beta\big):=\card\big(\M(g)\cap[\alpha,\beta]\big)
    -(\beta-\alpha)\card\big(\M(g)\big)\,,
  \end{equation*}
where $0\le\alpha\le\beta\le 1$.
In order to prove that $\M(g)$ is approximately uniformly distributed,
which is the same as saying that $\J$ is \emph{uniformly distributed} in
$[1,p]$, we have to show that the \emph{extreme discrepancy}
  \begin{equation*}
    \D\big(\M(g)\big):=\frac{1}{\card\big(\M(g)\big)}
    \sup_{1\le\alpha\le\beta\le 1}
    \big|\D\big(\M(g); \alpha,\beta\big)\big|    
  \end{equation*}
becomes small when $p$ gets large. This is the object of the following
theorem.
\begin{theorem}\label{Theorem1}
There exist absolute constants $c_1,c_2>0$, such that if 
$\frac 1\pi < p(\beta-\alpha)$, then 
  \begin{equation}\label{eqTh1}
    \big|\D\big(\M(g); \alpha,\beta\big)\big|\le c_1p^{1/2}\log p
    \big(2+\log p(\beta-\alpha)\big)
  \end{equation}
for any $0\le\alpha\le\beta\le 1$, and
  \begin{equation}\label{eqTh2}
    \big|\D\big(\M(g)\big)\big|\le
    \frac{c_2}{\card\big(\M(g)\big)}\cdot
    p^{1/2}\log^2 p\,.
  \end{equation}
\end{theorem}

A consequence of Theorem~\ref{Theorem1} assures us that any interval whose length 
combined with the length of $\J$ exceeds a certain margin, contains plenty of
elements of $\LL(g)$. 
 
%%%%%%%%%%%%%%%%%
\begin{corollary}\label{Corollary1}
For any $\delta>0$, any subinterval of $[0,p-1]$ of length $M$ contains
at least $(1-\delta)\frac{MN}{p}$ and at most $(1+\delta)\frac{MN}{p}$
elements of $\LL(g,\J)$, provided that $MN>\frac{c_3}{\delta}p^{3/2}\log^2 p$ 
for some absolute constant $c_3>0$.
\end{corollary}

%%%%%%%%%%%%%%%%%%%%%%%%%%%%%%%%%%%%%%%%%%%% 
\section{Estimate of an Exponential Sum}\label{sec:expsum}
One way to get bounds for the discrepancies is to obtain estimates for
certain exponential sums (see \eqref{eq7} below), and this our first point.

Let $\theta$ and $\zeta$ be roots of unity of order $p-1$ and $p$,
respectively. 
We consider the twisted sum, called the Lagrangian resolvent of $\theta$ and
$\zeta$:
  \begin{equation}\label{eq1}
    S(\theta,\zeta):=\zeta+\theta\zeta^g+\cdots+\theta^{p-2}\zeta^{g^{p-2}}\,.
%    =\sum_{j=0}^{p-2}\theta^j\zeta^{g^j}
  \end{equation}
Plainly $S(1,1)=p-1$ and it is known that 
  \begin{equation}\label{eq2}
    S(\theta,\zeta)\le \sqrt p\,,
  \end{equation}
for all $\theta$ and $\zeta$ that are not both equal to $1$. 
Let us see this for completeness. We have:
\begin{equation*}
\begin{split}
|S(\theta,\zeta)|^2=&
\sum_{k=0}^{p-2}\sum_{l=0}^{p-2}\theta^{k-l}\zeta^{g^{k}-g^{l}}\\
=& p-1+
\sum_{k=0}^{p-2}\sum_{\substack{l=0\\l\neq k}}^{p-2}
\theta^{k-l}\zeta^{ g^l(g^ { k-l } -1)  }\,.
\end{split}
\end{equation*}
Let us see that here, for any $l$ fixed, the differences $k-l$ run over the set of nonzero classes
$\mod (p-1)$. Then, since the order of both $\theta$ and $g$  is $p-1$, the sums above are equal to
\begin{equation*}
\begin{split}
=&
\sum_{t=1}^{p-2}\theta^t
	\sum_{l=0}^{p-2} \zeta^{ g^l(g^ { t } -1)  }
		=\sum_{t=1}^{p-2}\theta^t
		\sum_{s=1}^{p-1} \zeta^{ s  }\\
=&\sum_{t=1}^{p-2}\theta^t \cdot(-1) =1\,.
\end{split}
\end{equation*}
and \eqref{eq2} follows.

By \eqref{eq1}, we get
  \begin{equation}\label{eq3}
    \sum_{j=0}^{p-2}\theta^{kj}\zeta^{u(g^j-z)}=\zeta^{-uz}
    S(\theta^k,\zeta^u)\,.
  \end{equation}
Then we sum relations \eqref{eq3} over $1\le u\le p$. Note that
  \begin{equation*}
    \sum_{u=1}^{p}\zeta^{u(g^j-z)}=
    \begin{cases}
      p, &     \text{if $g^j\equiv z\pmod p$;}\\
      0, & \text{otherwise},
    \end{cases}
  \end{equation*}
and observe that since $0\le j\le p-2$, the condition 
$g^j\equiv z\pmod p$ can be written as $j=\log_g z$. These yield
  \begin{equation}\label{eq4}
    \theta^{k\log_g z}=\frac 1p 
    \sum_{u=1}^{p}\zeta^{-uz} S(\theta^k,\zeta^u)\,.    
  \end{equation}
Now taking $\theta=e_{p-1}(1)$ and $\zeta=e_p(1)$, where
$e_q(x):=\exp\big(\frac{2\pi i x}{q}\big)$, and summing equalities
\eqref{eq4} over $z\in\J$, we obtain
  \begin{equation}\label{eq5}
    \sum_{z\in\J}e_{p-1}\big(k\log_g z\big)=\frac 1p 
    \sum_{u=1}^{p}S(\theta^k,\zeta^u)
    \sum_{z\in\J}e_{p}(-uz)\,.    
  \end{equation}
The sum over $z$ on the right-hand side is sharply bounded by
   \begin{equation}\label{eq6}
     \begin{split}
       \bigg|\sum_{z\in\J}e_p(-uz)\bigg|
       &=\bigg|\sum_{j=1}^Ne_p\big(u(a+jr)\big)\bigg|
       =\bigg|\sum_{j=1}^Ne_p(ujr)\bigg|\\
       &\le \min \bigg( N,\ \frac{2}{\big|e_p(ur)-1 \big|} \bigg)
       \le \min \bigg( N,\ \frac{1}{\big|\sin  \frac{\pi u r}{p} \big|} \bigg)\\
       & \le \min \Big( N,\ \Big( 2\,\norm{\frac{ur}{p}}\Big)^{-1} \Big)\,,
\end{split}
\end{equation}
where $\norm{\cdot}$ is the distance to the nearest integer.
Then, using the \eqref{eq2} and \eqref{eq6} in \eqref{eq5}, we
conclude that
   \begin{equation}\label{eq7}
     \begin{split}
       \bigg|\sum_{z\in\J}e_{p-1}\big(k\log_g z\big)\bigg|
       &\le\frac 1p\sum_{u=1}^pp^{1/2}\cdot
       \min \Big( N,\ \Big( 2\,\norm{\frac{ur}{p}}\Big)^{-1} \Big)\\
       &\le\sqrt p + p^{-1/2}\sum_{u=1}^{p-1}
       \bigg(2\norm{\frac{ur}{p}}\bigg)^{-1}\\
       &\le\sqrt p + p^{-1/2}\sum_{v=1}^{\frac{p-1}2} \frac{p}{v}
       \le \sqrt p(2+\log p)\,.
     \end{split}
\end{equation}
The estimate \eqref{eq7} is slightly more general than the
P\'olya-Vinogradov inequality for character sums.

%%%%%%%%%%%%%%%%%%%%%%%%%%%%%%%%%%%%%%%%%%%%%%%
\section{The Proof of Theorem~\ref{Theorem1} and Corollary~\ref{Corollary1}}\label{sec:proof}
A bound for the discrepancy can be deduced applying the Erdo\"os-Tur\'an
inequality~\cite[Chapter 1, page 8]{Montgomerycarte}.  
This says that for any $0\le \alpha\le\beta\le 1$ and any positive
integer $K$, we have 
  \begin{equation*}
    \begin{split}
      \big|\D(\M(g);\,\alpha,\,\beta)\big|
      \leq& \frac{|\M(g)|}{K+1}+ 2\sum_{k=1}^{K}
                \bigg ( \frac{1}{K+1}+ 
                  \min \Big (\beta-\alpha,\ \frac{1}{\pi k} \Big )\bigg )
                  \bigg | \sum_{x\in\M(g)}\exp(2\pi ikx)\bigg|\,.
    \end{split}
  \end{equation*}
Bounding the exponential sum by \eqref{eq7}, the right-hand side is
  \begin{equation*}
    \begin{split}
      \leq& \frac{|\M(g)|}{K+1}+ 2\sqrt p(2+\log p)\bigg(1+\sum_{k=1}^{K}
                  \min \Big (\beta-\alpha,\ \frac{1}{\pi k} \Big )\bigg)\\
      \leq& \frac{|\M(g)|}{K+1}+ 2\sqrt p(2+\log p)\bigg(1+\sum_{1\le k\le \frac{1}{\pi(\beta-\alpha)}}^{}(\beta-\alpha)
                  +\sum_{\frac{1}{\pi(\beta-\alpha)}< k\le K}^{}\frac{1}{\pi k}\Big )\bigg)\\
      \leq& \frac{|\M(g)|}{K+1}+ c\sqrt p\log p\bigg(1+\Big|\log K(\beta-\alpha)\Big|\bigg)\,,
    \end{split}
  \end{equation*}
for some absolute constant $c>0$.
If we take $K=p-1$ in this estimate, we obtain \eqref{eqTh1}.

Next, let us see that if $\beta-\alpha\le 1/\pi p$, then $\M(g)$ contains at most one element,
therefore
  \begin{equation}\label{eq8}
    \begin{split}
      \frac{1}{\card(\M(g))}\big|\D(\M(g);\,\alpha,\,\beta)\big|
      \leq&\frac{1}{\card(\M(g))}\Big(1+(\beta-\alpha)\card(\M(g))\Big)\\
	\leq&\frac{2}{\card(\M(g))}\,.
    \end{split}
  \end{equation}
When $\beta-\alpha>1/\pi p$, we apply \eqref{eqTh1}, and obtain
  \begin{equation}\label{eq9}
    \begin{split}
      \frac{1}{\card(\M(g))}\big|\D(\M(g);\,\alpha,\,\beta)\big|
      \leq&\frac{1}{p}+\frac{c\sqrt p \log^2 p}{\card(\M(g))}\\
	\leq&\frac{c'\sqrt p\log^2 p}{\card(\M(g))}\,,
    \end{split}
  \end{equation}
 for some absolute constant $c'>0$.
Now \eqref{eqTh2} follows from \eqref{eq8} and \eqref{eq9}, and this concludes the proof of Theorem~\ref{Theorem1}.
\smallskip

To prove Corollary~\ref{Corollary1}, 
let $\I=[s,t]\subset [0,p-1]$ be any subinterval of length $t-s=M>0$, and let
$\delta>0$. 
Let $\alpha=s/p$ and $\beta=t/p$. 
We may assume that $\delta>1/\sqrt p$, 
since otherwise the result is trivial. Let $\alpha, \beta\in[0,1]$ with $\beta-\alpha=N/p$.

By the hypothesis, it follows that $\sqrt p\log^2 p<c''MN/p$ for some $c''>0$, 
and then by Theorem~\ref{Theorem1}, it implies that 
$\big|\D(\M(g);\,\alpha,\,\beta)\big|\le c''MN/p$. 
This can be rewritten as
  \begin{equation*}
    \begin{split}
	(1-\delta)\frac{MN}{p}\le \card(\M(g)\cap[\alpha,\beta])\le(1+\delta)\frac{MN}{p}\,,
     \end{split}
  \end{equation*}
which proves the corollary.

%%%%%%%%%%%%%%%%%%%%%%%%%%%%%%%%%%%%%%%%%%%%%%%
\section{A few open problems}\label{sec:problems}

There are different points of view and ways to study the distribution of the
elements of a certain sequence. But going further along the lines followed
above, let us first notice that Theorem~\ref{Theorem1} and
Corollary~\ref{Corollary1} applies not only to $\M(g,\J)$, but for sets
featuring certain patterns such as those generated when $\J$ is replaced by
unions of arithmetic progressions, also. This is easy to see, since 
\begin{equation*}
	\D(\M_1\cup\M_2)=\D(\M_1)+\D(\M_2)\,,
\end{equation*}
for any sets $\M_1,\M_2\subset[0,1]$ with $\M_1\cap\M_2=\emptyset$.

A further step in the evaluation of changes produced by the discrete logarithm
function would be to evaluate the discrepancy when the original set ($\J$ in the
notation from the introduction) is additionally changed by a non linear
transform. Such an example would require to estimate a sum such as, for
instance,
\begin{equation*}
\sum_{x\in\J} e_{p-1}\big(P(\log_g x)\big)\,,
\end{equation*}
where $P(x)=a_0+a_1x+\cdots+a_nx^n$, with $a_0,\dots,a_n$ integers,
$a_n\not\equiv 0\pmod p$ and $n\ge 2$.

Another spreading factor appears if more than one primitive root are
involved. Let $g_1,\dots,g_r$ be primitive roots $\mod p$ and let 
$a,b_1,\dots, b_r$ be integers. Then the problem is to find a nontrivial
estimate for the sum
\begin{equation*}
\sum_{x\in\J}
e_{p-1}\big(ax+b_1\log_{g_1}x+b_2\log_{g_2}x+\cdots+b_r\log_{g_r}x\big)\,.
\end{equation*}

Related to these questions is the problem that asks to study the changes
produced by the discrete logarithm function in the order of its arguments.
If the elements of $\LL(g,\J)$ were randomly distributed in $[0,p-1]$, then
comparing the size, for $x_1, x_2\in\J$ with $x_1<x_2$, one expects that about
half of the time $\log_{g}x_1<\log_g x_2$ and half of the time
$\log_{g}x_1>\log_g x_2$.  And similarly, for any fixed $r\ge 2$, when
$p\to\infty$,  all the $r!$ possible arrangements among the numbers 
$\log_g x_1,\dots,\log_g x_r\in [0,p-1]$ should occur with about the same
frequency when $(x_1,\dots,x_r)$ runs over $\J^r$.

%%%%%%%%%%%%%%%%%%%%%


\begin{thebibliography}{999}
\bibitem{DH} W. Diffie and M. Hellman,
\emph{New directions in cryptography},
IEEE Trans. Inform. Theory {\bf 22} (1976), 644--654.

\bibitem{ElGamal}
T. ElGamal,
\emph{A public key cryptosystem and a signature scheme based on
  discrete logarithms},
IEEE Trans. Inform. Theory {\bf 31} (1985), 469--472.

\bibitem{DSA}
D. W. Kravitz, 
\emph{Digital signature algorithm}, 
U.S. patent \#5,231,668, 27 July, 1993.

\bibitem{Menzes} A. Menzes, P. C. Van Oorschot and S. A. Vanstone,
\emph{Handbook of Applied Cryptography},
CRC Press, 1996, xxxvi+780 pp.

\bibitem{Montgomerycarte}
H. L. Montgomery, 
{\em Ten lectures on the interface between analytic
number theory and harmonic analysis,}
CBMS, {\bf 84.} Providence, RI, 1994, xiv+220 pp.

\bibitem{Odlyzko} A. Odlyzko,
\emph{Discrete Logarithms: The Past and the Future},
Preprint 1999.

\bibitem{BS}
W. D. Banks and I. E Shparlinski, \textit{Exponential sums with polynomial values of the discrete logarithm}, 
Uniform Distribution Theory, \textbf{2} (2007), no. \textbf{2}, 67–-72.

\end{thebibliography}
\end{document}